\documentclass[12pt]{article}
\usepackage{amssymb}
\usepackage{float}
\restylefloat{figure}

\newtheorem{theorem}{Theorem}[section]
\newtheorem{proposition}[theorem]{Proposition}

\newtheorem{definition}[theorem]{Definition}

\newcommand{\sign}{\mathrm{sign}}

\newcommand{\mod}{\mathrm{mod}}

\begin{document}
\title{ Monte Carlo simulations with a generalized detailed balance using the quantum-classical isomorphism }

\author{Yefim I. Leifman \\
Department of Mathematics and Statistics, \\ Bar-Ilan University, \\ 
52900, Ramat-Gan, Israel \\ leifmany@lycos.com }

\maketitle

\begin{abstract}

The main idea of this work is that the quantum-classical isomorphism is a suitable framework for a  generalization of the notion of detailed balance. The quantum-classical isomorphism is used in order to develop a Monte Carlo simulation with controlled deviation from detailed balance, that is with a generalized detailed balance and known relative entropy with respect to the reference process at each point. In order to apply this method to molecular simulations a new algorithm for realization of a partial chirotope, based on linear programming, a new distance geometry algorithm and a new all-atom off-lattice Monte Carlo method are proposed.  

{\bf Keywords:} Detailed balance; Distance geometry; Monte Carlo simulation;  Quantum-classical isomorphism; Relative entropy

\end{abstract}

\section{Background} \label{introduction} 
 
The characteristic time of an event in the molecular world is $10^{-15}$ sec, i.e. one iteration of a {\it molecular dynamics} must simulate changes corresponding to a time of this order. Biomolecules of interest, such as proteins, have thousands of atoms and, even using the simplest approximation for the molecular potential and a powerful computer, only thousand to million of iterations can be performed in a reasonable time. The approximate time of protein folding even in vivo is $10^{-5}$ sec. The timing for a standard Amber benchmark 159 residue protein in water is 249 ps/day of simulations on a single 3.4 GHz processor  \cite{amber-06}. The performance of the NAMD program on different platforms can be viewed in \cite{phillips-05}. However, significant changes of molecular conformations were achieved, for example, in all-atom molecular dynamics simulations of 36 residue protein on supercomputer of hundreds processors \cite{duan-98}. Another method of molecular simulation is molecular {\it Monte Carlo} simulation.  It does not try to simulate the physical movement of a molecule but only visit ({\it sample}) its conformational space according to an appropriate probability distribution, such as Boltzmann distribution. It gives a hope of closing the aforementioned time gap.

A Monte Carlo (MC) simulation is a random process. Usually, it is a {\it Markov process} \cite{ito-93}. The Markov process $\{X_t\}_{t\in T}$ is specified by its {\it transition probability} $P(s,x,t,Y)=P(X_t\in Y|X_s=x)$, $s\leq t$,  and {\it initial distribution} - the distribution of the random variable $X_0$. $Y$ is in the smallest $\sigma $-algebra that contains all open sets of the state space $S$ of the process. If the transition probability depends only on the difference between $s$ and $t$, that is, if there exists a function $P(t,x,Y)$, such that $P(s,x,t,Y)=P(t-s,x,Y)$, then the Markov process is called {\it temporally homogeneous}. A {\it jump process} is a continuous time process which changes its state after non-zero time. A temporally homogeneous Markov jump process is determined by its initial distribution, jump rate $j(x)$ for a  jump from each $x\in S$ and (jump) transition probability $P(x,Y)$ with (jump) transition probability density $p(x,y)$. Only temporally homogeneous Markov jump processes with $j(x)=1$ for all $x\in S$ will be considered in this work. In this case jumps occur according to the Poisson process. 

A process with the transition probability density $p(x,y)$ satisfies the {\it semi-detailed balance} condition \cite{frank-05} if there exists a density $\mu $ such that for all $y\in S$
$$\int_{x\in S} \mu (x)p(x,y)dx=\mu (y).$$ 
One says that such process with the initial density $\mu $ is in a {\it steady state}.

A process with the transition probability density $p(x,y)$ satisfies the {\it detailed balance} condition if there exists a density $\mu $ such that for all $x,y\in S$
$$\mu (x)p(x,y)=\mu (y)p(y,x).$$
One says that such process with the initial density $\mu $ is in {\it equilibrium}.

If a process with a finite number of states satisfies the detailed balance condition, then the process converges to a limiting distribution, which is unique \cite{mano-99}, \cite{nara-01}.
The simplest example of a process with detailed balance is a process with independent outcomes, that is $p(x,y)=p(z,y)$ and $\mu(x)=p(z,x)$ for all $x,y,z\in S$.

 A {\it Metropolis Monte Carlo} simulation \cite{metropolis-53} is an example of a process with detailed balance. Let given a system of particles in 3-dimensional space. Let $r$ and $n$ be vectors of coordinates of two states of the particle system  in the phase space. The {\it Boltzmann law} gives the ratio of densities to be in the state $r$ or $n$ in equilibrium
$$\frac{\mu (n)}{\mu (r)}=\exp \Bigl( \frac{-(E(n)-E(r))}{kT}\Bigr) ,$$
where $E(r)$ is the energy of the conformation of the particles, $T$ is a temperature, $k$ is the Boltzmann constant. The aim of the Metropolis 
MC sampling algorithm is to sample the phase  
space according to the distribution which satisfies the Boltzmann law.  First, one generates a conformation $r$ (present). Next, one generates a new conformation $n$ by adding a small random displacement to $r$. One must now decide whether the new conformation will be accepted or rejected. One wants to choose a transition probability density $p(r,n)$ such that the detailed balance condition
 $$\mu (r)p(r,n)=\mu (n)p(n,r)$$
 is satisfied. 
We denote the probability density to try a move from $r$ to $n$ by 
$l(r,n)$ and the probability of accepting a move from $r$ 
to $n$ by $\alpha(r,n)$. Assume that the transition probability density is given by 
$$p(r,n)=\alpha (r,n)l(r,n).$$ 
If $l$ is symmetric, i.e. $l(r,n)=l(n,r)$ ,
 then detailed balance implies 
$$\mu (r)\alpha(r,n)=\mu (n)\alpha(n,r),$$ 
 and therefore, 
$$\frac{\alpha(r,n)}{\alpha(n,r)}=\exp \Big( \frac{-(E(n)-E(r))}{kT}\Big) .$$ 
One of the possibilities to satisfy this condition is the choice of 
Metropolis et al. \cite{metropolis-53}
$$\alpha(r,n)=\min\{1,\exp \Big(\frac{-(E(n)-E(r))}{kT}\Big)\}.$$
 If $E(n)\leq E(r)$, 
then the move is accepted. If $E(n)>E(r)$, then we generate a random 
number $U$ from the uniform distribution in the interval $[0,1]$ and 
we accept the move if $U<\alpha(r,n)$ and reject otherwise. $l$ is not specified, except for the assumption that it is symmetric. 
This reflects freedom in the choice of moves.

It is obvious that an infinite number of states of the gas corresponds to a given macroscopic condition of the gas.  Through macroscopic measurements one should not be able to distinguish between two gases existing in different states (thus corresponding to two distinct representative points in phase space) but satisfying the same macroscopic conditions. Thus when one speaks of a gas under certain macroscopic conditions, one is in fact referring not to a single 
state, but to an infinite number of states. In other words, one refers not to a single system, but to a collection of systems, identical in composition and macroscopic condition but existing in different states. Such a collection of systems is called an {\it ensemble} (Chapter 4 of \cite{huang-87}), which is geometrically represented by a distribution of representative points in phase space, usually a continuous distribution. An ensemble is completely specified by this distribution. Metropolis MC samples an ensemble which is called {\it canonical ensemble}. It is appropriate to a system whose temperature is determined through contact with a heat reservoir (Chapter 7 of \cite{huang-87}). 

In {\it generalized-ensemble} simulations \cite{okamoto}, each state is weighted by a non-Boltzmann probability weight factor. This allows the simulation to escape from energy barriers and to sample much wider  space than by conventional methods. Monitoring the energy in a single simulation run in such ensembles, one can obtain also canonical ensemble averages as functions of temperature. One of the best-known generalized-ensemble methods is the {\it replica-exchange} MC \cite{okamoto},\cite{pitera}. The system for a replica-exchange MC consists of non-interacting copies, or replicas, of the original system in canonical ensemble at different temperatures. There is a one-to-one correspondence between replicas and temperatures. A simulation of replica-exchange MC is then realized by alternately performing the following two steps. 
\begin{enumerate}
\item each replica in the canonical ensemble at a fixed temperature is simulated simultaneously and independently for a certain number of Metropolis MC steps,
\item A random pair of replicas which are at neighboring temperatures $T_m$ and $T_{m+1}$ are exchanged moving the low-temperature conformation to the high-temperature simulation and vice versa. These replica swaps are accepted according to the Metropolis criterion with the acceptance probability 
$$\alpha_m=\min \{1,\exp\Big(\frac{E_{m+1}}{kT_{m+1}}+\frac{E_m}{kT_m}-\frac{E_{m+1}}{kT_m}-\frac{E_m}{kT_{m+1}}\Big)\}.$$ 
\end{enumerate}
The effect of replica exchange is to prevent low-temperature simulations from becoming trapped in local minima, because they are occasionally swapped to higher temperatures where they can escape these minima and move to other regions of phase space. Simultaneously, the low-temperature simulations are always being seeded with low-energy conformations produced by simulations at higher temperatures. The replica-exchange MC sampling is in detailed balance.

The notion of the {\it quantum-classical isomorphism}  (sometimes it is called just the {\it classical isomorphism}) originated in Chapter 10 of \cite{feynman-65}. The derivation can be found in Chapter 10 of \cite{allen-87}. 
Consider a quantum particle at temperature $T$, $\beta =1/(kT)$. The density operator is $\rho =\exp(-\beta H)$, where $H=-\frac{\hbar^2}{2m}\frac{d^2}{dx^2}+V(x)$. The density matrix is
$$\rho (r,r',\beta)=<r|\exp (-\beta H)|r'>$$
$$=<r|\exp (-\beta H/K)\dots \exp (-\beta H/K)\dots \exp (-\beta H/K)|r'>.$$
Inserting unity in the form $1=\int |r><r|dr$ of the integral of projectors $|r><r|$ over the volume of the system between each exponential gives
$$\rho (r,r',\beta)=\int <r|\exp (-\beta H/K)|r_2><r_2|\exp (-\beta H/K)|r_3>$$
$$\dots <r_{K-1}|\exp (-\beta H/K)|r_K><r_K|\exp (-\beta H/K)|r'>dr_2\dots dr_K$$
$$=\int \rho(r,r_2,\beta /K)\rho(r_2,r_3,\beta /K)\dots \rho(r_K,r',\beta /K)dr_2\dots dr_K.$$
If $\beta /K$ is sufficiently small, the following approximation is valid \cite{feynman-65}
$$\rho (r_a,r_b,\beta /K)\approx \rho _{free}(r_a,r_b,\beta /K)\exp \Bigl(-\frac {\beta}{2K}(V^{cl}(r_a)+V^{cl}(r_b))\Bigr) $$
where $V^{cl}(r)$ is the classical potential energy, and the free-particle density matrix for a single particle of mass $m$ is \cite{feynman-65}
$$\rho _{free}(r_a,r_b,\beta /K)=\Bigl( \frac{mK}{2\pi \beta \hbar ^2}\Bigr) ^{\frac{D}{2}}\exp \Bigl(-\frac{mK}{2\beta \hbar ^2}r_{ab}^2\Bigr) ,$$
where $r_{ab}=|r_a-r_b|$ and $D$ is a dimension of the space. Therefore
$$\rho (r_1,r_1,\beta )\approx \Bigl( \frac {mK}{2\pi \beta \hbar^2}\Bigr) ^{\frac{DK}{2}}\int \exp \Bigl( -\frac {mK}{2\beta \hbar ^2}(r_{12}^2+r_{23}^2+\dots r_{K1}^2)\Bigr)$$
\begin{equation} \label{density}
 \exp \Bigl( -\frac{\beta}{K}(V^{cl}(r_1)+V^{cl}(r_2)+\dots V^{cl}(r_K))\Bigr) dr_2 \dots dr_K.
\end{equation}
The quantity $\rho (r,r,\beta)$ is proportional to the density to be in $r$. If an $N$-particle system is considered, replace $D$ by $ND$. That is, the density of the $N$-particle quantum system which satisfies the Boltzmann law corresponds to the density of the $KN$-particle classical system which satisfies the Boltzmann law and consists of $K$ copies of the $N$-particle classical systems and  in the copies with neighbor numbers $i$ and $i+1$ for $1\leq i<K$ and $K$ and $1$ the corresponding particles  are connected by springs with spring constant $mK/(\beta^2 \hbar ^2)$. This approximation becomes exact as $K\rightarrow \infty $ \cite{allen-87} and can be used in a conventional MC simulation to investigate quantum properties.

The main difficulties in all-atom detailed balanced MC simulations of biochemical processes  involving big molecules  are that these methods have huge autocorrelation time (Section 2 of \cite{sokal-96}) and these processes, generally, do not approach to equilibrium.

A process is {\it reversible} if  $p(x,y)>0$ implies $p(y,x)>0$ for all $x,y\in S$.
Let us mention some characteristics of non-equilibrium reversible processes.
Such process satisfies the {\it master equation}
$$\frac{d\mu _t(x)}{dt}=\int_{y\in S}\mu _t(y)p(y,x)dy-\mu _t(x).$$
For the probability density $\mu _t$ the {\it Gibbs entropy} is given by
$$S_G(\mu _t)=-\int_{x\in S}\mu _t(x)\log \mu _t(x)dx.$$
$dS_G(\mu _t)/dt$ is deduced from master equation and according to Section 2.4 of \cite{lebowitz-99} it splits  into the {\it entropy production rate} $R(\mu _t)$ and the {\it entropy flow rate} $A(\mu _t)$ as follows
$$\frac{dS_G(\mu _t)}{dt}=R(\mu _t)-A(\mu _t),$$
where
$$R(\mu _t)=\frac{1}{2}\int \int_{x,y\in S}(\mu _t(x)p(x,y)-\mu _t(y)p(y,x))\log \frac{\mu _t(x)p(x,y)}{\mu _t(y)p(y,x)}dxdy$$
and
$$A(\mu _t)=\int_{x\in S}\mu _t(x)I(x)dx=\langle I \rangle _{\mu _t},
 I(x)=\int_{y\in S}p(x,y)\log \frac{p(x,y)}{p(y,x)}dy.$$
The entropy production rate is expressed in Section 3.1 of \cite{gaspard-05} in terms of the "particle fluxes"
 $$J_t(x)=\mu _t(x)p(x,y)-\mu _t(y)p(y,x)$$
 and "forces"
 $$F_t(x)=\log \frac{\mu _t(x)p(x,y)}{\mu _t(y)p(y,x)}.$$
As was mentioned in \cite{blr-05}, entropy production rate can be considered as a measure of a lack of equilibrium. 

\begin{definition} [\cite{gray-90}] \label{relative}
For two probability distributions $P$ and $Q$ with probability densities $p$ and $q$, the {\it relative entropy} (Kullback-Leibler divergence) is defined by 
$$D_{KL}(P\| Q)=\int_{x\in S}p(x)\log \frac{p(x)}{q(x)}dx.$$
If entropy is measured in bits, the logarithm in this formula is taken to base 2, or to base $e$, if entropy is measured in nats.
\end{definition}

The Gibbs inequality says that $D_{KL}(P\| Q)\geq 0$ and the relative entropy is zero iff $P=Q$. The entropy flow rate $A(\mu _t)$ for the distribution $\mu _t$ which is concentrated at one point $x$ is the relative entropy of the transition probability density $p'(y)=p(x,y)$ and the probability density $p''(y)=p(y,x)$.

\section{The problems addressed in this work}  \label{introduction2}

In this work we consider three different, but related problems:
\begin{enumerate}
\item generalizing the detailed balance condition in order to include processes which do not preserve any distribution, but with the property that, roughly speaking, it is known how much information the "jumper" must retrieve from its path in order to reach its current position (when the corresponding reference process is given), then building processes with such generalized detailed balance,
\item enhancing all-atom off-lattice molecular MC simulations which are in detailed balance, for example replica-exchange MC, in order that they can be performed with a move set consisting of separate moves of each atom, with all degrees of freedom also in the case of dense atom packing; we do this by means of a new distance geometry algorithm, which plays in such MC simulations a role which is similar to the role of the SHAKE algorithm \cite{ryckaert-77}, \cite{todorov-06}  in molecular dynamics,
\item building an initial sample for all-atom off-lattice molecular MC simulations according to chirality constraints and distance constraints, and additionally, geometric manipulations with a molecule, which preserve, as far as possible, the aforementioned constraints, but exploit flexibility of a molecule.
\end{enumerate}

The connecting link of the following considerations of these problems is the distance geometry procedure which we call "centering". We build an initial sample of a molecule satisfying molecular chirality constraints and distance constraints with the help of linear programming and subsequent "iterative vibrant centering". We perform a Metropolis MC with a move set consisting of separate moves of each atom in the sample space which is restricted with the help of "centering". We combine this   distance geometry method and this MC to exploit flexibility of a molecule in geometric manipulations with it. We use the same restricted sample space in MC without detailed balance.

Section \ref{generalized} provides an example of a process with generalized detailed balance with respect to its reference process using the quantum-classical isomorphism.

Section \ref{numerical} provides the results of the numerical experiments considering some properties of the process which are similar to the example from Section \ref{generalized}, but can be numerically examined.

Clearly, there is a temperature at which the Amber force field 
cannot ensure the integrity of the molecule in the all-atom Metropolis MC sampling with a move set consisting of separate moves of each atom as described in \cite{metropolis-53}. This limits replica temperatures in replica-exchange MC. The same phenomenon can happen in non-equilibrium simulations. A way to overcome this difficulty is to change the potential in order to restrict such deformations. We change the potential with the help of distance geometry "centering" procedure. It can be considered as building the restricted sample space which includes all relevant conformations. It is described in Section \ref{restricted}.

In order to start a Metropolis MC with this potential one has to build an initial sample which is in the restricted sample space. As will be proved in Section \ref{iterative}, if one starts at some point of the sample space and performs as many centerings as needed from an infinite sequence of centerings  which contains an infinite number of centerings of each atom, then at certain  step one reaches a point in the aforementioned restricted sample space. We call this algorithm "iterative centering".

Molecular chirality constraints impose limitations on molecular conformations. These limitations are in addition to the limitations imposed by the weighted graph of the desired distances. A new algorithm for  realization of a partial chirotope, based on linear programming is proposed in Section \ref{realization}.

Suppose, that a sample which satisfies a given partial chirotope (chirality constraints) is built. Now we need to push  this sample into the restricted sample space. Numerical tests show, that reiteration of iterative centering  algorithm with chirality checking can become jammed if the initial sample is far from the restricted sample space. For overcoming this difficulty the  "vibrant iterative centering" distance geometry algorithm is proposed in Section \ref{metropolis}. The vibrant iterative centering algorithm can be incorporated in a convenient computation scheme with a Metropolis MC in the restricted sample space. This scheme allows flexible manipulations with a molecule with further equilibration. It is described in Section \ref{metropolis}.

Section \ref{molecular} describes a molecular MC simulation without detailed balance using the quantum-classical isomorphism. In this simulation we use a process which is similar to the example  which is described in Section \ref{generalized}. The simulation of Section \ref{molecular} is only  preliminary since the choice of appropriate parameters for true molecular simulations is not considered in this work, but the observations of Section \ref{molecular} hints at the possibility to use such method in molecular MC.

All algorithms with centerings are new. The iterative vibrant centering can be useful in existing distance geometry software in order to improve its sampling properties. The restricted sample space can be useful in all-atom off-lattice molecular simulation software, for example, in order to increase the temperature of the hottest replica in replica-exchange MC. The algorithm for realization of a partial chirotope using linear programming is also new. It can be useful with iterative vibrant centering. 

The notion of generalized detailed balance in a framework of Langevin dynamics was proposed in \cite{kurchan-97}. As far as we know, our work is the first attempt to generalize the detailed balance condition in order to include processes which are not in a steady state but with the property, that it is known how much information the jumper must retrieve from its path in order to reach its current position (when the corresponding reference process is given). In our opinion  the appropriate generalization of the detailed balance condition is the most important problem in computer molecular simulations, since we believe that an all-atom detailed balanced simulation of working ribosome will never be possible on a digital computer. 

The concepts of vertex, particle, bead and atom represent the same object in different contexts of this work. The notations of Sections  \ref{restricted}, \ref{iterative}, \ref{realization}, \ref{metropolis} are independent of the notations of Sections \ref{introduction},  \ref{generalized}, \ref{numerical},  \ref{molecular}.

\section{Generalized detailed balance} \label{generalized}

Consider the Young's double-slit experiment (\S\S 26,27,32,33 of \cite{matveev}). Let two parallel plane screens be separated by a distance $l$. Let there be $N\geq 1$ rectilinear parallel slits of width $s$ on the first screen. In the case of the Young's double-slit experiment $N=2$. Let the correspondent slit borders of neighbor slits be separated by a distance $d$. Let there be a plane monochromatic light wave with length $\lambda$ which propagates perpendicularly to the screens and hits the second screen through the slits in the first screen.  
If $\lambda \ll s\leq d\ll l$, then the intensity on the second screen is  
$$I(y)=\frac{I_0\sin^2Y\sin^2(NqY)}{N^2Y^2\sin^2(qY)},$$
where $y$ is a distance between a point on the second screen and perpendicular projections of the slits on the second screen,
$$Y=\frac{\pi s}{\lambda l}y,$$ 
$q=d/s\geq 1$ and 
$I_0$ is the maximum intensity on the second screen. 
In this case the Fraunhofer diffraction takes place, that is $I(ly_1)/I(ly_2)$ does not depend on $l$. 

Since 
$$\int_{0}^{\infty}\frac{\sin^2(ax)}{x^2}dx=|a|\frac{\pi}{2}$$  (3.821.9 from \cite{gradshteyn}), without loss of generality let the intensity distribution for $N=1$ be 
$$I_1(x)=\frac{\sin^2x}{\pi x^2}.$$
Since 
$$\int_{0}^{\infty}\frac{\sin^2(ax)\cos^2(bx)}{x^2}dx=\frac{a\pi}{4}$$ for $0<a\leq b$ (3.828.11 from \cite{gradshteyn}), let the intensity distribution for $N=2$ be 
$$I_2(x)=\frac{\sin^2x\sin^2(2qx)}{2\pi x^2\sin^2(qx)}.$$

Then $$\int_{-\infty}^{\infty}I_1(x)\log_2\frac{I_1(x)}{I_2(x)}dx=1$$
since
$$\int_{0}^{\infty}\frac{\log \cos^2(ax)}{x^2}\cos (bx)dx=-a\pi+\pi b\log 2+\pi\sum_{n=1}^{m}\frac{(-1)^n(b-2an)}{n},$$
where $m\leq b/(2a)<m+1$, $m=0,1,2,3,\dots$
(I.158 from \cite{oberhettinger}).

In this example 1 bit of the information through what slit the photon passed translates in 1 bit of relative entropy $D_{KL}(I_1(x)\| I_2(x))$ for all $q=d/s\geq 1$. The process "without information" can be considered as a reference process. The formula for Fraunhofer diffraction is an approximation, but this property can be directly verified in experiments on various screens. We call this property by {\it the balance of relative entropy}.

In the path integral formulation of quantum mechanics \cite{feynman-65} the contribution of a particular path in the total probability amplitude for a photon has a phase proportional to time to travel along this path. In the aforementioned example the contribution of a particular path changes along this path such that the total intensity obeys the property of balanced relative entropy. A time with such property "rotates" in order to hide the information about the past which was not retrieved in time.

The rest of this section is devoted to giving an example of the process with jumps which obey the property of balanced relative entropy. If the reference process is in detailed balance, then the condition of the balance of relative entropy can be considered as a generalization of the detailed balance condition.

In order that the quantum-classical isomorphism be valid, the Boltzmann law must be satisfied \cite{feynman-65}.   Therefore MC with independent outcomes of one quantum particle using the quantum-classical isomorphism with condition, that the vector of the coordinates $a_0$ of the first copy of the particle is in the zero point must be defined by the following recurrent  {\it Levy construction} \cite{levy-39} for $1\leq n <K$
\begin{equation} \label{levy}
a_n=\frac{K-n}{K-n+1}a_{n-1}+\frac{1}{\sqrt{\frac{mK}{\beta \hbar ^2}(1+\frac{1}{K-n})}}\eta_n
\end{equation}
where each $\eta_n$ is a vector with independent standard normal distributed coordinates.
The Levy construction samples intermediate time points of a Brownian motion, conditioned to arrive to a predetermined point after a predetermined time.
A {\it Brownian bridge} is a Brownian motion conditioned to return to the initial point after a predetermined time. The aforementioned MC sample $(a_0,\dots ,a_{K-1})$ is the Levy construction for a Brownian bridge  for the time interval divided into $K$ equal subintervals with $a_0=0$. In the case of the quantum-classical isomorphism, this Brownian bridge time parameter is considered as "imaginary time" \cite{feynman-65} in contrast with ordinary time where the jumps of MC take place. In this work, imaginary time is discrete and denoted by a subscript  and ordinary time is continuous, but a number of jumps in jump processes is denoted by a superscript.

Consider MC simulations of two distinguishable quantum particles. Firstly we define an auxiliary  process $\{N_t\}$.

Let $(a_0,\dots ,a_{K-1})$ and $(b_0,\dots ,b_{K-1})$ be the aforementioned Levy constructions with $a_0=0$ and $b_0=0$. Shift them by $a$ and $b$ according to the distribution of $N_0$, that is $(a_0^1,\dots ,a_{K-1}^1)=(a_0,\dots ,a_{K-1})+(a,\dots ,a)$ and $(b_0^1,\dots ,b_{K-1}^1)=(b_0,\dots ,b_{K-1})+(b,\dots ,b)$. This is an initial sample for $\{N_t\}$. 
  
Suppose, that $i$ samples of $\{N_t\}$ are built, that is $i-1$ jumps of $\{N_t\}$ have already happened. Build new Levy constructions $(\hat a_0^i,\dots ,\hat a_{K-1}^i)$ and $(\hat b_0^i,\dots ,\hat b_{K-1}^i)$. Randomly choose two numbers $n_1^i$ and $n_2^i$, $n_1^i<n_2^i$, from 0 to $K-1$ with equal probability for each pair. If 
\begin{equation} \label{lessormore}
\| a_{n_1^i}^i-b_{n_1^i}^i\| < \| a_{n_2^i}^i-b_{n_2^i}^i\| 
\end{equation}
shift these new Levy constructions $(\hat a_0^i,\dots ,\hat a_{K-1}^i)$ and $(\hat b_0^i,\dots ,\hat b_{K-1}^i)$ by $a_{n_1^i}^i-\hat a_{n_1^i}^i$ and $b_{n_2^i}^i-\hat b_{n_2^i}^i$ correspondingly with probability $\alpha$, otherwise, shift them by $a_{n_2^i}^i-\hat a_{n_2^i}^i$ and $b_{n_1^i}^i-\hat b_{n_1^i}^i$. If $\| a_{n_1^i}^i-b_{n_1^i}^i\| > \| a_{n_2^i}^i-b_{n_2^i}^i\| $ shift them by $a_{n_2^i}^i-\hat a_{n_2^i}^i$ and $b_{n_1^i}^i-\hat b_{n_1^i}^i$ with probability $\alpha$, otherwise, shift them by $a_{n_1^i}^i-\hat a_{n_1^i}^i$ and $b_{n_2^i}^i-\hat b_{n_2^i}^i$ and so on. 

$\{N_t\}$ satisfies (\ref{density}) if it is in a steady state.  
Denote $\{N_t\}$ with $\alpha=1$ by $\{M_t\}$ and denote $\{N_t\}$ with $\alpha=1/2$ by $\{L_t\}$. The process $\{L_t\}$ is a well-known MC of two free particles using quantum-classical isomorphism.   Jumps of $\{M_t\}$ satisfies the balanced relative entropy condition with respect to jumps of $\{L_t\}$. Now, given a process with aforementioned property and a detailed balanced process used as a reference process, we know how much information a jumper must retrieve from its path in order to reach its current position.

\section{Numerical experiments} \label{numerical}

In order to present the results of the numerical experiments, we define some processes which are similar to $\{N_t\}$. All these processes consider $K$ copies of two particles in $\mathbb{R}$. The copies  are connected by springs as described in the discussion of the quantum-classical isomorphism in Section \ref{introduction}. The first sample for all the  processes of this section is $((0,\dots ,0),(\delta,\dots,\delta))$, where $\delta>0$ is  large relative to the lengths of the steps of the processes. 

Define  $\{N_t^j\}$ for $0<j<K$ as follows. Suppose, that $i$ samples of $\{N_t^j\}$ are built. Build new Levy constructions $(\hat a_0^i,\dots ,\hat a_{K-1}^i)$ and $(\hat b_0^i,\dots ,\hat b_{K-1}^i)$ according to (\ref{levy}). Randomly choose a number $n_1^i$  from 0 to $K-1$ with equal probability for each number. Let $n_2^i\equiv n_1^i+j \; (\mod \;K)$. If 
$\| a_{n_1^i}^i-b_{n_1^i}^i\| < \| a_{n_2^i}^i-b_{n_2^i}^i\| $
shift these Levy constructions $(\hat a_0^i,\dots ,\hat a_{K-1}^i)$ and $(\hat b_0^i,\dots ,\hat b_{K-1}^i)$ by $a_{n_1^i}^i-\hat a_{n_1^i}^i$ and $b_{n_2^i}^i-\hat b_{n_2^i}^i$ correspondingly with probability $\alpha$ (these jumps we call "forward jumps") otherwise shift them by $a_{n_2^i}^i-\hat a_{n_2^i}^i$ and $b_{n_1^i}^i-\hat b_{n_1^i}^i$ ("backward jumps"). If $\| a_{n_1^i}^i-b_{n_1^i}^i\| > \| a_{n_2^i}^i-b_{n_2^i}^i\| $ shift them by $a_{n_2^i}^i-\hat a_{n_2^i}^i$ and $b_{n_1^i}^i-\hat b_{n_1^i}^i$ with probability $\alpha$ ("forward") otherwise shift them by $a_{n_1^i}^i-\hat a_{n_1^i}^i$ and $b_{n_2^i}^i-\hat b_{n_2^i}^i$ ("backward") and so on. 

Define  $\{W_t^j\}$ for $0<j<K$ as follows. Suppose, that $i$ samples of $\{W_t^j\}$ are built. Randomly choose a number $n_1^i$  from 0 to $K-1$ with equal probability for each number. Let $n_2^i\equiv n_1^i+j \; (\mod \;K)$ and $h_q=mK/(\beta ^2\hbar ^2)$.

 If $\| a_{n_1^i}^i-b_{n_1^i}^i\| < \| a_{n_2^i}^i-b_{n_2^i}^i\| $, then with probability $\alpha$
$$a_{n_2^i}^{i+1}=\frac{a_{n_2^i-1}^i+a_{n_2^i+1}^i}{2}+\frac{1}{\sqrt{2\beta h_q}}\eta_i,\quad
b_{n_1^i}^{i+1}=\frac{b_{n_1^i-1}^i+b_{n_1^i+1}^i}{2}+\frac{1}{\sqrt{2\beta h_q}}\eta'_i$$
and other coordinates unchanged, otherwise
$$a_{n_1^i}^{i+1}=\frac{a_{n_1^i-1}^i+a_{n_1^i+1}^i}{2}+\frac{1}{\sqrt{2\beta h_q}}\eta_i,\quad
b_{n_2^i}^{i+1}=\frac{b_{n_2^i-1}^i+b_{n_2^i+1}^i}{2}+\frac{1}{\sqrt{2\beta h_q}}\eta'_i$$
and other coordinates unchanged.

 If $\| a_{n_1^i}^i-b_{n_1^i}^i\| > \| a_{n_2^i}^i-b_{n_2^i}^i\| $, then with probability $\alpha$
$$a_{n_1^i}^{i+1}=\frac{a_{n_1^i-1}^i+a_{n_1^i+1}^i}{2}+\frac{1}{\sqrt{2\beta h_q}}\eta_i,\quad
b_{n_2^i}^{i+1}=\frac{b_{n_2^i-1}^i+b_{n_2^i+1}^i}{2}+\frac{1}{\sqrt{2\beta h_q}}\eta'_i$$
and other coordinates unchanged, otherwise
$$a_{n_2^i}^{i+1}=\frac{a_{n_2^i-1}^i+a_{n_2^i+1}^i}{2}+\frac{1}{\sqrt{2\beta h_q}}\eta_i,\quad
b_{n_1^i}^{i+1}=\frac{b_{n_1^i-1}^i+b_{n_1^i+1}^i}{2}+\frac{1}{\sqrt{2\beta h_q}}\eta'_i$$
and other coordinates unchanged.

The following numerical experiments show that in some sense $\{W_t^j\}$ and $\{N_t^j\}$ are similar.
Consider $K=8,16,32,64,128$, $\alpha =1,2/3,7/12,13/24,25/48$. Let $m$ be atomic mass unit,  $T=300$K. Standard normal distributed random numbers are obtained by the Box-Muller algorithm \cite{devroye-86} from the uniformly distributed pseudo-random numbers \cite{matsumoto-98}. We compute the following quantities: the number of performed jumps $\cal J$, the number $\cal F$ of forward jumps $i\leq \cal J$  for which 
$\| a_{n_1^i}^i-b_{n_1^i}^i\| - \| a_{n_2^i}^i-b_{n_2^i}^i\| $ and 
$\| a_{n_1^i}^{i+1}-b_{n_1^i}^{i+1}\| - \| a_{n_2^i}^{i+1}-b_{n_2^i}^{i+1}\| $ have different signs,
the number of such backward jumps $\cal R$, the average of the coordinates $ {\cal A}=\frac{1}{K}\sum_{n=1}^{K}a_n^{{\cal J}+1}$ after the last jump, the average ${\cal B}=\frac{1}{K}\sum_{n=1}^{K}(b_n^{{\cal J}+1}-\delta)$, the average number ${\cal C}=\frac{1}{K}\sum_{n=1}^{K}C_n$, where $C_n$ is the number of $i\leq \cal J$ such that 
$\| a_n^i-b_n^i\| - \| a_{n'}^i-b_{n'}^i\| $ and
$\| a_n^{i+1}-b_n^{i+1}\| - \| a_{n'}^{i+1}-b_{n'}^{i+1}\| $ have different signs, where $n'\equiv n+j \; (\mod \;K)$, $\frac{({\cal A}+{\cal B})K\sqrt{K}}{{\cal F}-{\cal R}}$ and  $\frac{({\cal A}+{\cal B}){\cal J}}{({\cal F}-{\cal R}){\cal C}}$.

The results for $\{N_t^j\}$ are shown in Table \ref{table1} of Appendix, the results for $\{W_t^j\}$ are shown in Table \ref{table2} of Appendix. In the case of $\{N_t^j\}$ we take ${\cal J}/{\cal C}=2$. The lengths in these tables are in units $\frac{ \hbar}{\sqrt{mkT}}$. The observations are as follows.
\begin{enumerate}
\item 
For both $\{N_t^j\}$ and $\{W_t^j\}$, if we fix $K$ and $j$, then $\cal A$ and $\cal B$ are proportional to $\alpha-1/2$.
\item  For $\{W_t^j\}$, if we fix $K$ , then  $\frac{{\cal A}+{\cal B}}{{\cal F}-{\cal R}}$ is approximately constant. It is true also when $j$ is not a constant, but a random variable, for example, uniformly distributed in some interval, as for $\{N_t\}$ (Not shown in the Appendix). 
\item 
If we fix $j/K$, the closeness of $\frac{({\cal A}+{\cal B}){\cal J}}{({\cal F}-{\cal R}){\cal C}}$ for $\{N_t^j\}$ and  for $\{W_t^j\}$ justifies that ${\cal J}/{\cal C}$ jumps of $\{W_t^j\}$  approximate two jumps of $\{N_t^j\}$.
\end{enumerate}

\section{The restricted sample space} \label{restricted}

Evidently, there is a temperature at which the Amber force field 
cannot ensure the integrity of the considered molecule in the all-atom Metropolis MC sampling with the move set consisting of the separate moves of each atom as described in \cite{metropolis-53}. This limits replica temperatures in the replica-exchange MC. The same phenomenon can happen in non-equilibrium simulations. A way to overcome this difficulty is to change the potential in order to restrain such deformations. 

A finite undirected weighted graph $G$ is a triple $<V,E,W>$, where $V$ 
denotes the set of its vertices, $E$ denotes the set of its edges, 
and $W:E \rightarrow {\mathbb R}^+$ is a function which specifies 
a positive weight for each graph edge.
In order to restrict the sample space as was mentioned in Section \ref{introduction2}, one has  to set a weighted graph which corresponds to the molecule. Generally, atoms are vertices of this graph, covalent bonds form a part of its edges, pairs of atoms which are bonded through two covalent bonds form another part of its edges, weights are desired distances. The weights of edges which connect two atoms  which are bonded through two covalent bonds determine the bond angles. So, the weighted graph of methane has 10 undirected weighted edges. The weighted graph of amide plane has also $C^\alpha -C^\alpha$ and $H-O$ edges, since their distances are well defined in amide plane. If one knows additional distances between atoms (for example from Nuclear Magnetic Resonance data), one adds corresponding edges and weights too.

Let $f:V \rightarrow {\mathbb R}^k$ be a {\it conformation} of $G$ in the $k$-dimensional Euclidean space ${\mathbb R}^k$ (regardless of the weights). Let the coordinates of vertices $a_v=f(v)$ be all distinct and let $h:E \rightarrow {\mathbb R}^+$ be the spring constants of edges. The point
\begin{equation} \label{center}
c_u=\frac{1}{\sum_{\{ v| \{u,v\}\in E\}}h(\{u,v\})}\sum_{\{ v|\{u,v\}\in E\} }h(\{u,v\})
\Big( a_v+\frac{W(\{u,v\})}{\| a_u-a_v\| } (a_u-a_v)\Big)
\end{equation} 
will be called the {\it center of the vertex (atom)} $u$.

The corresponding algorithm for finding the center of the vertex $u$ takes as its input the adjacency-list representation \cite{cormen-90} of the finite weighted graph $G=<V,E,W>$ and current coordinates of its vertices. The adjacency-list representation of the graph consists of the array $Adj$ of $|V|$ lists, one for each vertex in $V$. For each $u\in V$ the adjacency list $Adj[u]$ contains all the vertices $v\in V$ such that there exists an edge $\{u,v\} \in E$. The weight $W(\{ u,v\} )$ of the edge $\{ u,v\} \in E$ is stored with a vertex $v$ in $u$'s adjacency list. $h(\{ u,v\} )$ are stored like $W(\{ u,v\} )$ and the vector $A[u]$ of current coordinates of the vertex $u$ is stored with $u$. If $x$ is a pointer to an element of the list $Adj[u]$, then, according to pseudocode conventions \cite{cormen-90}, $vertex[x]$ denotes a vertex which adjacent to $u$ ( denote it by $v$ ), $W[x]$ and $H[x]$ denote $W(\{ u,v\} )$ and $h(\{ u,v\} )$.

\begin{tabbing}
\ \ \ \ \ \ \ \ \ \= nthen \= nthen \= nthen \= nthen \= nthen \= nthen \kill
\> $Center(u)$ \\
1\> $t \leftarrow 0$ \\
2\> $q \leftarrow 0$ \\
3\> $y \leftarrow (0,0,0)$ \\
4\> $f \leftarrow 0$ \\
5\> $x \leftarrow head(Adj[u])$ \\

6\> while $x \neq $ NIL \\
7\> \ do \> $v \leftarrow vertex[x]$ \\
8\>  \> $z \leftarrow A[u]-A[v]$ \\
9\>  \> $r \leftarrow \| z \| $ \\
10\>  \> if $r>0$ \\
11\>  \> \ then \> $y \leftarrow y+H[x](A[v]+(W[x]/r)z)$ \\
12\>  \>  \> $t \leftarrow t+H[x]$ \\
13\>  \> \ else \> $f \leftarrow f+H[x]W[x]$ \\
14\>  \>  \> $q \leftarrow q+H[x]$ \\
15\>  \> $x \leftarrow next[x]$ \\

16\> if $q>0$ and $t>0$ \\
17\> \ then \> $z \leftarrow (1/t)y-A[u]$ \\
18\>  \> $r \leftarrow \| z \| $ \\
19\>  \> if $r>0$ \\
20\>  \> \ then \> $y \leftarrow y+qA[u]+(f/r)z$ \\
21\>  \>  \> $y \leftarrow (t/(t+q))y$ \\

22\> if $t>0$ \\
23\> \ then \> return $(1/t)y$ \\
24\> \ else \> return $A[u]$ \\
\end{tabbing}

As will be proved in Section \ref{iterative}, if one starts at arbitrary point of a sample space and performs as many centerings as needed from a sequence of centerings which contains an infinite number of centerings of each vertex, then at some step one achieves a point of a sample space such that for each $u\in V$ it holds that $\| a_u-c_u \| <S(u)$ for a given $S(u)>0$. Denote the set of such points by $\cal D(S)$. As will be viewed in Section \ref{iterative} the {\it Hooke potential} $\sum_{e\in E}\frac{h(e)}{2}(\| e\| -W(e))^2$ is large outside $\cal D(S)$. We change a given potential, for example the Amber force field, by assuming it infinitely large outside of $\cal D(S)$. $S(u)$ has to be not too small if one does not want to restrict a considerable part of molecular degrees of freedom.

\section{The iterative centering algorithm} \label{iterative}

{\it Distance geometry} is a part of computational geometry
which is devoted to the study of the existence or non-existence of an 
embedding satisfying the condition in the following definition 
as well as methods for construction of such embedding.  

\begin{definition} \label{d5.1}
Let $G=<V,E,W>$ be a finite undirected weighted graph, where $V$ 
denotes the set of its vertices, $E$ denotes the set of its edges, 
and $W:E \rightarrow {\mathbb R}^+$ is a function which specifies 
a positive weight for each graph edge. An embedding of $G$ in the 
$k$-dimensional Euclidean space ${\mathbb R}^k$ is a function 
$f:V \rightarrow {\mathbb R}^k$ such that for each edge $e=\{ v,w\} 
\in E$ one has $\|f(v)-f(w)\|=W(e)$. $G$ is called $k$-embeddable 
iff such an embedding exists.
\end{definition}

The problem of $k$-Embeddability of an integer-weighted 
undirected graph is NP-hard \cite{saxe-79}. However there is 
a semidefinite programming algorithm for the Euclidean distance 
matrix completion problem, i.e. determining whether there 
exists a number $k$ for which a given undirected weighted 
graph is $k$-embeddable \cite{laurent-00}. If there exists an embedding according to Definition \ref{d5.1}, it is evidently in $\cal D(S)$.

\begin{theorem} [\cite{blumenthal-53}] \label{t5.2}
A complete graph is $k$-embeddable iff each of its complete subgraphs
 with $k+3$ vertices is k-embeddable.
\end{theorem}

Another distance geometry problem is that of {\it bounded} 
$k$-{\it Embeddability}, namely whether for given bounds 
$l,u:E\rightarrow {\mathbb R}^+$ there exists a weight 
$W$ with $l(e)\leq W(e)\leq u(e)$ for which a graph 
$G=<V,E,W>$ is $k$-embeddable.

There are number of methods which can be applied to solving
 the aforementioned problems, which generally arise in 
Nuclear Magnetic Resonance data interpretation: {\it metric matrix} 
distance geometry \cite{crippen-88},\cite{nmrchitect-98}, {\it simulated annealing}, 
{\it variable target} function optimization \cite{braun-85} and {\it global continuation}
\cite{more-97}.

The {\it iterative centering} distance geometry algorithm consists of performing as many steps as needed from an infinite sequence of centerings which contains an infinite number of centerings of each vertex.

\begin{proposition} \label{t5.3}
Let $\{ u_k \} $ be an infinite sequence of vertices of $G$. Let us apply a sequence $\{ Center(u_k) \} _{1\leq k\leq n}$ of $n$ centerings. Denote the displacement of the vertex $u_n$ after the iteration $Center(u_n)$ by $l_n$. Then $\lim_{n \to \infty}\| l_n\| =0$.
\end{proposition}
{\bf Proof. } Without loss of generality suppose that all coordinates of $u_n$  after $n-1$ iterations are zero. Denote the coordinate vectors of vertices adjacent to $u_n$ after $n-1$ iterations by $a_1,...,a_m$, the corresponding weights by $w_1,...,w_m$ and the spring constants by $h_1,...,h_m> 0$. If $\| a_i \| >0$  for every $i$,$1\leq i\leq m$, then according to (\ref{center}) the coordinates of $u_n$ after $n$ iterations is 
$$y=\frac{1}{\sum_{i=1}^{m}h_i}\sum_{i=1}^{m}h_i\Big( 1-\frac{w_i}{\| a_i\| }\Big) a_i.$$ 
Therefore $$\Big( \sum_{i=1}^{m}h_i \Big) \| y\| ^2=\sum_{i=1}^{m}h_i\Big\| \Big( 1-\frac{w_i}{\| a_i\| } \Big) a_i\Big\|^2- \sum_{i=1}^{m}h_i\Big\| \Big( 1-\frac{w_i}{\| a_i\| }\Big) a_i-y\Big\|^2$$ similar to the Huygens theorem about momenta. 
$$\Big\| \Big( 1-\frac{w_i}{\| a_i\| }\Big)a_i-y\Big\| \geq \Big| \|a_i-y\|-w_i\Big| $$ by the triangle inequality and  $$\Big\| \Big( 1-\frac{w_i}{\| a_i\| } \Big) a_i \Big\|=\Big| \| a_i\| -w_i \Big| .$$
  Therefore 
\begin{equation} \label{e5.3.1}
\Big( \sum_{i=1}^{m}h_i \Big) \| y\| ^2\leq \sum_{i=1}^{m}h_i (\| a_i\|-w_i)^2-\sum_{i=1}^{m}h_i(\|a_i-y\| -w_i)^2.
\end{equation}

If $\| a_i\| >0$ for every $i$, $1\leq i\leq k$, and $\| a_i\| =0$ for every $i$, $k+1\leq i\leq m$, we take some $z$, $\| z\| =1$ and we define 
$$y=\frac{1}{\sum_{i=1}^{m}h_i} \Big( \sum_{i=1}^{k}h_i \Big ( 1-\frac{w_i}{\| a_i\| } \Big) a_i+\sum_{i=k+1}^{m}h_iw_iz\Big).$$
Similarly to the previous, we have (\ref{e5.3.1}). We put 
$$z=\frac{1}{\Big\| \sum_{i=1}^{k}h_i \Big ( 1-\frac{w_i}{\| a_i\|}\Big) a_i  \Big\| }\sum_{i=1}^{k}h_i \Big ( 1-\frac{w_i}{\| a_i\|}\Big) a_i$$ as in the $Center(u)$ algorithm.
The right hand side of (\ref{e5.3.1}) is twice the difference of the old and the new values of the Hooke potential $\sum_{e\in E}\frac{h(e)}{2}(\| e\| -W(e))^2$. That is, at stage $n$, the Hooke potential decreases by at least $\Big( \sum_{i=1}^{m}h_i \Big) \| l_n\| ^2$, where $l_n$ is the shift of the center. If $\| l_n\| \not\rightarrow 0$, then the Hooke potential would become negative at some $n$. Since Hooke potential is non-negative we have $\lim_{n \to \infty}\| l_n\| =0$.    $\square$

\section{The realization of the partial chirotope  
related \\ to the molecular chirality constraints} \label{realization}

Let $x_1,x_2,\dots,x_r,y_1,y_2,\dots,y_r\in {\mathbb R}^r$. 
Then the following Grassmann-Plucker relation holds \\
$\det(x_1,x_2,\dots,x_r)\cdot \det(y_1,y_2,\dots,y_r)$
$$=\sum_{i=1}^r\det(y_i,x_2,\dots,x_r)
\cdot \det(y_1,\dots,y_{i-1},x_1,y_{i+1},\dots,y_r).$$
The difference of the left and the right sides is an alternating multi-linear 
form in the $r+1$ arguments $x_1,y_1,y_2,\dots,y_r$, which are 
vectors in an $r$-dimensional vector space; hence, the difference
 of the left and the right sides is identically zero. 
For example, in rank 3 one gets for every set of 5 vectors 
(denoting determinants by square brackets, and labeling 
the points 1 to 5) the relation $[123][145]-[124][135]+[125][134]=0$. 
This requires that these 6 signs of the brackets on the left side are 
such that the equality is at least possible for this sign pattern, 
when actual scalars are not given: for example, these 6 signs 
could be +,+,+,+,+,+, but not +,+,-,+,+,+.

\begin{definition} [\cite{bjorner-93}] \label{d6.1}
Let $r\geq 1$ be an integer, and let $E$ be a finite set. 
A chirotope of rank $r$ on $E$ is a mapping 
$\chi :E^r \rightarrow\{-1,0,1\}$  which satisfies the following 3 properties:
\begin{enumerate}
\item  $\chi$ is not identically 0,
\item  $\chi$ is alternating, that is, $\chi (a_{\sigma_1},a_{\sigma_2},\dots,a_{\sigma_r})=\sign(\sigma)\chi (a_1,a_2,\dots,a_r)$ for all $a_1,a_2,\dots,a_r\in  E$ and every permutation  $\sigma$,
\item for all  $a_1,a_2,\dots,a_r,b_1,b_2,\dots,b_r\in E$ such that \\
$\chi(a_1,a_2,\dots,a_r)\cdot\chi(b_1,b_2,\dots,b_r)\neq 0$, \\
there exists an $i\in\{1,2,\dots,r\}$ such that \\
$\chi (b_i,a_2,\dots,a_r)\cdot \chi(b_1,\dots,b_{i-1},a_1,b_{i+1},\dots,b_r)$ \\
$=\chi(a_1,a_2,\dots,a_r)\cdot\chi(b_1,b_2,\dots,b_r)$.
\end{enumerate}
\end{definition}

The axioms comes from abstracting sign properties in the 
Grassmann-Plucker relations for $r$-order determinants. The 
chirotope axioms are a version of the {\it oriented matroid} axioms. 

Suppose $E=\{1,2,\dots,n\}$. Given any $(n-r)$-tuple $(a_1,\dots,a_{n-r})$ 
of elements in $E$, then we write $(a_1',\dots,a_r')$ for some permutation 
of $E\backslash\{a_1,\dots,a_{n-r}\}$. Then $(a_1,\dots,a_{n-r},a_1',\dots,a_r')$ 
is a permutation of $(1,2,\dots,n)$, and we can compute \\ 
$\sign(a_1,\dots,a_{n-r},a_1',\dots,a_r')$ as the parity of the number of 
inversions of this string. The mapping $\chi ^*: E^{n-r} \rightarrow \{-1,0,1\}$, 
defined by  $$\chi ^*(a_1,\dots,a_{n-r}) = \chi (a_1',\dots,a_r')\sign(a_1,
\dots,a_{n-r},a_1',\dots,a_r')$$ is called the chirotope {\it dual} to the chirotope $\chi$.

Let $\chi : E^r\rightarrow \{-1,0,1\}$ be a chirotope, $E=\{1,\dots ,n\}$. 
If there exists $\{x_1,\dots,x_{n}\}\subset {\mathbb R}^r$ such that 
$$\chi(a_1,a_2,\dots,a_r)=\sign(\det(x_{a_1},x_{a_2},\dots,x_{a_r}))$$ 
for all $1\leq a_1<a_2<\dots<a_r\leq n$, then $\chi$ is called {\it realizable} 
and $q:E\rightarrow {\mathbb R}^r$, $i\mapsto x_i$ is called a 
{\it realization} of $\chi$.

Let $G_r({\mathbb R}^n)$ be the real 
Grassmann manifold of $r$-dimensional linear subspaces in 
${\mathbb R}^n$, or equivalently $Mat_{r\times n}({\mathbb R})/GL_r({\mathbb R})$, which corresponds to the space of configurations of $n$ vectors in ${\mathbb R}^r$ modulo the action of the general linear group $GL_r({\mathbb R})$. Thus the realization $q$  of $\chi$ corresponds to a point in $G_r({\mathbb R}^n)$. The set of such points is called the {\it realization space} of $\chi$.

Let $\chi : E^r\rightarrow \{-1,0,1\}$ be a chirotope. 
Then for each subset $\{a_1,\dots,a_{r+2}\}$ of $E$ there 
exists $\{x_1,\dots,x_{r+2}\}\subset {\mathbb R}^r$ such that 
$$\chi(a_{i_1},a_{i_2},\dots,a_{i_r})=\sign(\det(x_{i_1},x_{i_2},
\dots,x_{i_r}))$$ for all $1\leq i_1<i_2<\dots<i_r\leq r+2$. 
This feature of chirotopes is called by {\it local realizability}. 
 Local realizability 
follows from the facts that realizability is preserved under 
duality since $G_r({\mathbb R}^n)=G_{n-r}({\mathbb R}^n)$ and that all rank 2 chirotopes are realizable. 

The realizability problem for chirotopes is NP-hard \cite{bjorner-93}.
There is an algorithm for a realization of a chirotope 
$\chi$ which is valid when the realization 
space of $\chi$ is contractible and $\chi :E^r \rightarrow\{-1,1\}$ \cite{bokowski-86}, \cite{crippen-88}. 

If the alternating map $\chi $ is only partially defined and the Grassmann-Plucker relation holds whenever $\chi $ is defined on all its  participants, then $\chi $ is called a {\it partial chirotope}.
A partial chirotope $\chi '$ of rank $r$ on $E$ is called {\it extendable} if there exists a chirotope $\chi $ of rank $r$ on $E$ and for any $a_1,\dots ,a_r\in E$, $\chi' (a_1,\dots ,a_r)=\chi (a_1,\dots ,a_r)$ holds whenever $\chi '(a_1,\dots ,a_r)$ is defined.
The problem of testing extendability of a partial chirotope is NP-complete \cite{tschirschnitz}. 

Molecular chirality constraints impose limitations on molecular conformations.  These limitations are in addition to the limitations imposed by weighted graph of desired distances. 
A set of inequalities of type $\det(x_b-x_a,x_c-x_a,x_d-x_a)>0$, where $a,b,c,d\in V$, corresponds to molecular chirality constraints. Then the corresponding equalities $\chi (a,b,c,d)=1$ define a rank 4 partial chirotope. If $a\mapsto (x_a^1,x_a^2,x_a^3)$ satisfies these inequalities, then $a\mapsto (1,x_a^1,x_a^2,x_a^3)$ is a realization of the corresponding partial chirotope.
The most widely adopted method to realize a partial chirotope related to 
molecular chirality constraints is the minimization of the function, 
which includes deviations from given oriented volumes, by simulated 
annealing starting from an approximate embedding \cite{nmrchitect-98}. 
An example of a realization of a partial chirotope by use of such 
function can be found in \cite{vendruscolo-95}. 

Let "maximize $cx$ with conditions $Ax\leq b$ and $x\geq 0$" be a linear program. $x\geq 0$ means $x_j\geq 0$ for all $j$. Particularly, let $c_j$ be the price per unit of the product $j$ produced, $x_j$ be the quantity of the product $j$ produced, $b_i$ be the quantity of the material $i$ on hand, $a_{ij}$ be the  quantity of the material $i$ required to produce one unit of the product $j$. Let $y_i$ be the price per unit of the material $i$. One is interested in selling the materials instead of the products if $A^Ty\geq c$. The dual linear program "minimize $by$ with conditions $A^Ty\geq c$ and $y\geq 0$" answers the question what is the minimal price of all materials when it is advantageous to sell the materials instead of to work. This price is the same as the maximal income in the first (primal) linear program. It is the figurative formulation of the linear programming strong duality theorem as economists learn it.

In some cases the following algorithm allows one to realize a given molecular  partial chirotope.  Let $Y\subset V$ be a set of vertices, whose coordinates appear in inequalities of type $\det(x_b-x_a,x_c-x_a,x_d-x_a)>0$. Without loss of generality one can demand $\det(x_b-x_a,x_c-x_a,x_d-x_a)\geq \epsilon >0$ for all these inequalities and $z_i=x_i^3\geq 0$ for all $i\in Y$. If we fix $x_i^1$ and $x_i^2$ for all $i\in Y$ then the inequalities become linear. A feasible solution of the following (symmetric) linear programming problem 
$$\mbox{Minimize } z_1+...+z_m \mbox{ subject to } Az\geq \epsilon \mbox{ and } z\geq 0$$
is a solution of our problem.
Its dual problem
$$\mbox{Maximize } \epsilon t_1+...+\epsilon t_k \mbox{ subject to } A^Tt\leq 1 \mbox{ and } t\geq 0$$ 
has zero as a feasible solution. If the original problem has a feasible solution, then its dual is bounded by the strong duality theorem. The dual problem can be solved by the primal simplex procedure and if it is bounded, then the solution of the original problem can be taken from the last simplex tableau, according to the Chapter 4 of \cite{jeter}.

In this method one has to fix $x_i^1$ and $x_i^2$ for all $i\in Y$. Since $\epsilon$ is a minimal volume for a parallelepiped spanned by  $x_b-x_a,x_c-x_a,x_d-x_a$  for each ordered base $(a,b,c,d)$, the points $(x_a^1,x_a^2)$, $(x_b^1,x_b^2)$, $(x_c^1,x_c^2)$, $(x_d^1,x_d^2)$ cannot be on the same straight line. We place all atoms of a molecule in a sequence and  choose $x_i^1=\cos(2\pi i/n)$,  $x_i^2=\sin(2\pi i/n)$ for all $i\in V$ and $\epsilon =(\sin(2\pi /n))^3$, where $n$ is a number of atoms.

In practical implementation of this algorithm of realization of a molecular partial chirotope, one has to set a partial chirotope of a given molecule. For example, for $C^\alpha $ atom of amino acid residue it is necessary to set 3 ordered bases. Fixing only 2 of them jams vibrant iterative centering algorithm (the modification of iterative centering which will be described in Section \ref{metropolis}) and fixing 4 of them is too restrictive for the choice of $x_i^1$ and $x_i^2$. The fourth ordered base will be recovered by means of distance constraints. Similarly, for $C^\alpha $ atoms of one spire of $\alpha $-helix, in which participate 5 residues, it is necessary to set 3 ordered bases. This partial chirotope will be used also for chirality checking $CheckChirality(u)$. 

Consider an example of poly-L-threonine Thr$_{180}$. Each Thr residue contains two chiral centers. For $C^\beta $ atom of Thr residue it is necessary to set 3 ordered bases. Arrange 14 atoms of each Thr residue in a following order H-N-H-C$^\alpha$-C$^\beta$-H-O$^\gamma$-H-C$^\gamma$-H-H-H-C-O (or in the notations of Protein Data Bank H-N-H-CA-CB-H-OG1-H-CG2-H-H-H-C-O). The proposed algorithm successfully finds a realization of a corresponding partial chirotope. Let us add to this chirotope also the constraints on $C^\alpha $ atoms which appear assuming Thr$_{180}$ is twisted in 50 spires of right $\alpha $-helix. The algorithm successfully finds a realization in this case.

\section{Metropolis Monte Carlo in the restricted \\ 
sample space} \label{metropolis}

If a given partial chirotope is realized, one has to transform this realization, keeping correct chiralities, in order to achieve $\cal D(S)$ and then to start the Metropolis MC simulation in the restricted sample space. Let $CheckChirality(u)$ be a function which checks whether quadruples of vertices which contain a vertex $u$ satisfy a given partial chirotope. Numerical tests show that reiteration of the iterative centering  algorithm with chirality checking can become jammed if the initial sample is far from $\cal D(S)$. (For example, consider the weighted graph with four vertices and four edges on a plane: let the starting configuration be $A=(0,0)$, $B=(4,3)$, $C=(4,-3)$, $D=(0,5)$, the weights of $(A,B)$ and $(A,C)$ be 5, the weight of $(B,C)$ be 6, the weight of $(A,D)$ be 0.01, $B,C,D$ be counter-clockwise and $S(u)=0.001$ for all vertices.) For overcoming this difficulty one can use the following modification of the iterative centering algorithm:

\begin{tabbing}
\ \ \ \ \ \ \= nthen \= nthen \= nthen \= nthen \= nthen \= nthen \kill
\> $VibrantCenter(u)$ \\
1\> $a \leftarrow A[u]$ \\
2\> $z \leftarrow Center(u)$ \\
3\> $r \leftarrow \| a-z \| $ \\
4\> if $r>C\cdot S[u] $ \\
5\> \ then \> $A[u] \leftarrow a+C \cdot S[u] \cdot ((z-a)/r+c \cdot RandomVector())$ \\
6\> \ else \> if $r>S[u] $ \\
7\>  \> \ then \> $A[u] \leftarrow a+S[u] \cdot ((z-a)/r+c \cdot RandomVector())$ \\
8\>  \> \ else \> $A[u] \leftarrow z+S[u] \cdot c \cdot RandomVector()$ \\
9\> if not $CheckChirality(u)$ \\
10\> \ then \> $A[u] \leftarrow a$ \\
\end{tabbing}

$RandomVector()$ denotes a function, which returns a uniformly distributed vector in a sphere of radius $1$ with the center at the origin. A coefficient $c>1$ is introduced for ergodicity. A coefficient $C>1$ is introduced for speeding-up. There is no guarantee that $\cal D(S)$ will be achieved, but the examples of this section show that the method works.

The vibrant iterative centering algorithm can be useful if the method described in Section \ref{realization} fails. Split some vertex from $Y$ into several vertices and spread the  inequalities of the form $\det(x_b-x_a,x_c-x_a,x_d-x_a)\geq \epsilon $ in which the coordinates of the initial vertex participate over these new vertices.  In the weighted graph put the desired distances between the new vertices be 0. Do this for several vertices from $Y$. Apply the described linear programming method and the vibrant iterative centering to bring nearer the vertices obtained from the same vertex.

Introduce $CheckDistance$ function:
\begin{tabbing}
\ \ \ \ \ \ \= nthen \= nthen \= nthen \= nthen \= nthen \= nthen \kill
\> $CheckDistance(u)$ \\
1\> $x \leftarrow head(Adj[u])$ \\
2\> while $x \neq $ NIL \\
3\> \ do \> $v \leftarrow vertex[x]$ \\
4\>  \> if $\| A[v]-Center(v))\| <S[v]$ or $S[v]=0$\\
5\>  \> \ then \> $x \leftarrow next[x]$ \\
6\>  \> \ else \> return FALSE \\
7\> return TRUE \\
\end{tabbing}

Let us sum up the proposed methods in the following computing scheme:

\begin{tabbing}
\ \ \ \ \ \ \= nthen \= nthen \= nthen \= nthen \= nthen \= nthen \= nthen \= nthen \kill
\> $TrialMove(u)$ \\
1\> $a \leftarrow A[u]$ \\
2\> $A[u] \leftarrow a+S[u]\cdot RandomVector() $ \\
3\> if $\| A[u]-Center(u)\| <S[u]$ and $CheckDistance(u)$ and $CheckChirality(u)$\\
4\> \ then \> $E \leftarrow potential(u)$ \\
5\> \> $z \leftarrow A[u]$ \\
6\> \> $A[u] \leftarrow a$ \\
7\> \> $e \leftarrow potential(u)$ \\
8\> \> if $E<e$ or $Random()<\exp ((e-E)/(kT))$ \\
9\> \> \ then \> $A[u] \leftarrow z$ \\
10\> \ else \> $A[u] \leftarrow a$ \\
11\> \>  $z \leftarrow Center(u)$ \\
12\> \>  $r \leftarrow \| a-z \| $ \\
13\> \>  if $r>S[u]$ or not $CheckDistance(u)$ \\
14\> \> \ then \> if $r>C\cdot S[u] $ \\
15\> \> \> \ then \> $A[u] \leftarrow a+C \cdot S[u] \cdot ((z-a)/r+c \cdot RandomVector())$ \\
16\> \> \> \ else \> if $r>S[u] $ \\
17\> \> \> \> \ then \> $A[u] \leftarrow a+S[u] \cdot ((z-a)/r+c \cdot RandomVector())$ \\
18\> \> \> \> \ else \> $A[u] \leftarrow z+S[u] \cdot c \cdot RandomVector()$ \\
19\> \> \> if not $CheckChirality(u)$ \\
20\> \> \> \ then \> $A[u] \leftarrow a$ \\
\end{tabbing}

$Random()$ generates a uniformly distributed in $[0,1]$ random number \cite{matsumoto-98}. 

Now we apply the proposed methods in their natural succession: firstly a realization of partial chirotope, then vibrant iterative centering and then the Metropolis MC in the restricted sample space. If a conformation which satisfies chirality constraints is not given, then one can use the method of Section \ref{realization} to build such conformation. It is useful to rescale this conformation to its natural scale and proportions. Then one can start vibrant iterative centering with large $S[u]$ to break frozen parts of the initial conformation and gradually decrease $S[u]$ to required values. Then it is possible to start $TrialMove$ over all atoms, which drives a molecule to $\cal D(S)$ and then becomes to be the Metropolis MC simulation in $\cal D(S)$. Similar to the original Metropolis MC \cite{metropolis-53}, only coordinates of one atom are changed near its current position in a trial move. It makes the all-atom (including hydrogen) Metropolis MC possible also in the case of dense atom packing. Also $TrialMove$ allows flexible manipulations with molecule by adding or removing weighted edges of the weighted graph and subsequent equilibration.

Consider an example of poly-L-alanine Ala$_{36}$, which is twisted in 10 spires of right $\alpha $-helix. We add to the weighted graph, which is derived from the primary structure of the molecule, distances of hydrogen bonds $O_{(i)}-H_{(i+4)}$, $O_{(i)}-N_{(i+4)}$  and distances $C_{(i)}^\alpha -C_{(i+2)}^\alpha$, which are well defined in $\alpha $-helix (the parentheses contain numbers of residues). Then we apply the simplex procedure, vibrant iterative centering algorithm and Metropolis MC simulation in the restricted sample space using Amber force field as described in the previous sections and receive the expected structure.

\section[Molecular MC simulation without detailed 
balance \\ using the quantum-classical  isomorphism]{Molecular MC simulation without detailed \\
balance using the quantum-classical \\ isomorphism} \label{molecular}

In order to proceed to non-equilibrium molecular simulations we add pairs of particles similar to that described in the example of Section \ref{generalized} to the considered molecule. These pairs of particles are used as artificial devices and do not represent physical particles. In this section we shall call these artificial added particles by beads for convenience. Also we add one Hooke term per bead to the molecule potential so that it connects a bead to some atom of a molecule by a spring with spring constant $h_a$ and zero length when the spring is relaxed. Suppose, that a sample which satisfies a given partial chirotope and distance constraints of such equipped molecule is built and equilibrated by methods described in previous sections. Subsequently we produce $K$ copies of this system which are connected by springs as described in the discussion of the quantum-classical isomorphism in Section \ref{introduction}.

In this loaded case we cannot use two Levy constructions for an added pair of beads as described in Section \ref{generalized}, but in order to approach to a process whose jumps obey the property of the balance of  relative entropy let $K=2j$ and choose the $l$-th added pair of beads with probability $p_l$, then uniformly choose two copies  $n_1^i$ and $n_2^i=n_1^i+j$, fix the $n_1^i$-th copy of one bead from the $l$-th added pair and the $n_2^i$-th copy of another bead from the $l$-th added pair according to (\ref{lessormore}) and perform predefined large number of Metropolis MC steps over the rest of atoms and beads as described in Section \ref{metropolis} and so on.

Consider the example of the linear polymer molecule which contains $N=16$ identical atoms with some Lennard-Jones constants and with neighbor atoms connected by springs. We add one aforementioned pair of beads to each pair of neighbor atoms. We constrain the sequence of the second beads of the added pairs to have chiralities of right helix and produce two copies of this system which are connected by springs as described in Section \ref{introduction}. Then we start the simulation and observe that the polymer moves ahead. If one fixes the last atom in the space, then the polymer twists around the fixed atom like boa. The twisting polymer squeezes itself out. This observation hints on the possibility to use such method in molecular MC.

\newpage

\section*{Appendix} \label{appendix}

The following notations are used in Table \ref{table1} and Table \ref{table2} (see Section \ref{numerical} for details).

\begin{tabular}{lp{0.85\textwidth}}
$K$ & the number of copies in the quantum-classical isomorphism,\\
$j$ & the superscript parameter of $\{N_t^j\}$ or $\{W_t^j\}$,\\
$\alpha$ & the probability parameter of $\{N_t^j\}$ or $\{W_t^j\},$ \\
$\cal J$ & the number of forward and backward jumps, \\
$\cal F$ & the number  of forward jumps $i\leq \cal J$ for which 
$\| a_{n_1^i}^i-b_{n_1^i}^i\| - \| a_{n_2^i}^i-b_{n_2^i}^i\| $ and 
 $\| a_{n_1^i}^{i+1}-b_{n_1^i}^{i+1}\| - \| a_{n_2^i}^{i+1}-b_{n_2^i}^{i+1}\| $ have different signs, \\
$\cal R$ & the number  of backward jumps $i\leq \cal J$ for which 
$\| a_{n_1^i}^i-b_{n_1^i}^i\| - \| a_{n_2^i}^i-b_{n_2^i}^i\| $ and 
 $\| a_{n_1^i}^{i+1}-b_{n_1^i}^{i+1}\| - \| a_{n_2^i}^{i+1}-b_{n_2^i}^{i+1}\| $ have different signs, \\
$\cal A$ & the average of coordinates ${\cal A}=\frac{1}{K}\sum_{n=1}^{K}a_n^{{\cal J}+1}$ after the last jump, \\
$\cal B$ & the average ${\cal B}=\frac{1}{K}\sum_{n=1}^{K}(b_n^{{\cal J}+1}-\delta)$, where the first sample for the processes  is $((0,\dots ,0),(\delta,\dots,\delta))$,\\
$\cal C$ & the average number ${\cal C}=\frac{1}{K}\sum_{n=1}^{K}C_n$, where $C_n$ is the number of $i\leq \cal J$  such that $\| a_n^i-b_n^i\| - \| a_{n'}^i-b_{n'}^i\| $ and 
$\| a_n^{i+1}-b_n^{i+1}\| - \| a_{n'}^{i+1}-b_{n'}^{i+1}\| $  have different signs, where $n'\equiv n+j \; (\mod \;K)$, \\
\end{tabular}

In the case of $\{N_t^j\}$ we take ${\cal J}/{\cal C}=2$. The lengths in these tables are in units $\frac{ \hbar}{\sqrt{mkT}}$.

\newpage
\begin{table}[H] 
\begin{scriptsize}
\begin{tabular}{rrrrrrrrr}

K & j & $\alpha$ & ${\cal J}$ & ${\cal F}$ & ${\cal R}$ & ${\cal A}$ & ${\cal B}$ & $\frac{({\cal A}+{\cal B}){\cal J}}{({\cal F}-{\cal R}){\cal C}}$ \\
\\
8 & 2 & 1.0000 & 10000000 & 5002391 & 0 & 1.2212e+006 & 1.2242e+006 &  9.778e-001 \\ 
8 & 2 & 0.6667 & 10005000 & 3336249 & 1668173 & 4.0719e+005 & 4.0778e+005 &  9.773e-001 \\ 
8 & 2 & 0.5833 & 10000000 & 2919143 & 2081936 & 2.0472e+005 & 2.0419e+005 &  9.768e-001 \\ 
8 & 2 & 0.5417 & 20000000 & 5415950 & 4583516 & 2.0225e+005 & 2.0394e+005 &  9.758e-001 \\ 
8 & 2 & 0.5208 & 20000000 & 5205748 & 4792723 & 1.0111e+005 & 9.7748e+004 &  9.629e-001 \\ 

16 & 1 & 0.6667 & 10000000 & 3331214 & 1667291 & 2.2680e+005 & 2.2837e+005 &  5.471e-001 \\ 
16 & 2 & 0.6667 & 10000000 & 3333718 & 1666593 & 3.1050e+005 & 3.1231e+005 &  7.470e-001 \\ 
16 & 4 & 0.6667 & 10000000 & 3332375 & 1668024 & 4.0726e+005 & 4.0487e+005 &  9.763e-001 \\ 
16 & 8 & 0.6667 & 10000000 & 3332983 & 1667959 & 4.7057e+005 & 4.6881e+005 &  1.129e+000 \\ 

32 & 1 & 0.6667 & 10000000 & 3333992 & 1665310 & 1.6566e+005 & 1.6342e+005 &  3.946e-001 \\ 
32 & 2 & 0.6667 & 10000000 & 3334183 & 1666927 & 2.2823e+005 & 2.2691e+005 &  5.461e-001 \\ 
32 & 4 & 0.6667 & 10000000 & 3335290 & 1664622 & 3.1246e+005 & 3.1090e+005 &  7.460e-001 \\ 
32 & 8 & 0.6667 & 10000000 & 3334115 & 1665497 & 4.0766e+005 & 4.0893e+005 &  9.788e-001 \\ 

64 & 2 & 0.6667 & 10000000 & 3333389 & 1665903 & 1.6270e+005 & 1.6609e+005 &  3.944e-001 \\ 
64 & 8 & 0.6667 & 10000000 & 3334884 & 1667409 & 3.1174e+005 & 3.1077e+005 &  7.465e-001 \\ 
64 & 16 & 0.6667 & 10000000 & 3333361 & 1665149 & 4.0878e+005 & 4.0681e+005 &  9.783e-001 \\ 

\\
\end{tabular}
\end{scriptsize}
\caption{Numerical results for $\{N_t^j\}$} 
\label{table1}
\end{table}

\newpage
\begin{table}[H]
\begin{tiny}
\begin{tabular}{rrrrrrrrrrr}

K & j & $\alpha$ & ${\cal J}$ & ${\cal F}$ & ${\cal R}$ & ${\cal A}$ & ${\cal B}$ & ${\cal C}$ & $\frac{({\cal A}+{\cal B})K\sqrt{K}}{{\cal F}-{\cal R}}$ & $\frac{({\cal A}+{\cal B}){\cal J}}{({\cal F}-{\cal R}){\cal C}}$ \\
\\
16 & 1 & 0.6667 & 10000000 & 2294496 & 1153377 & 1.5090e+004 & 1.5224e+004 & 512272 & 1.700e+000 & 5.187e-001 \\ 
16 & 2 & 0.6667 & 10000000 & 1642125 & 827546 & 1.1095e+004 & 1.1140e+004 & 369874 & 1.747e+000 & 7.381e-001 \\ 
16 & 4 & 0.6667 & 10000000 & 1243440 & 631437 & 8.4986e+003 & 8.5709e+003 & 282097 & 1.785e+000 & 9.887e-001 \\ 

32 & 1 & 0.6667 & 10000000 & 2251225 & 1130260 & 5.2767e+003 & 5.2874e+003 & 251470 & 1.706e+000 & 3.748e-001 \\ 
32 & 2 & 0.6667 & 10000000 & 1580763 & 798450 & 3.7948e+003 & 3.7869e+003 & 178242 & 1.754e+000 & 5.436e-001 \\ 
32 & 3 & 0.6667 & 10000000 & 1307517 & 663010 & 3.1831e+003 & 3.1756e+003 & 148053 & 1.786e+000 & 6.665e-001 \\ 
32 & 4 & 0.6667 & 10000000 & 1138660 & 582221 & 2.7568e+003 & 2.7539e+003 & 129549 & 1.792e+000 & 7.644e-001 \\ 
32 & 5 & 0.6667 & 10000000 & 1042861 & 534477 & 2.5251e+003 & 2.5549e+003 & 118876 & 1.809e+000 & 8.405e-001 \\ 
32 & 6 & 0.6667 & 10000000 & 967434 & 495375 & 2.3921e+003 & 2.3670e+003 & 110380 & 1.825e+000 & 9.131e-001 \\ 
32 & 7 & 0.6667 & 10000000 & 911386 & 467753 & 2.2203e+003 & 2.2288e+003 & 104017 & 1.815e+000 & 9.644e-001 \\ 
32 & 8 & 0.6667 & 10000000 & 874623 & 449650 & 2.1260e+003 & 2.1291e+003 & 100035 & 1.812e+000 & 1.001e+000 \\ 
32 & 9 & 0.6667 & 12000000 & 1024333 & 527608 & 2.5445e+003 & 2.5039e+003 & 117076 & 1.840e+000 & 1.041e+000 \\ 
32 & 10 & 0.6667 & 11000000 & 911136 & 470663 & 2.2247e+003 & 2.2062e+003 & 104318 & 1.821e+000 & 1.061e+000 \\ 
32 & 11 & 0.6667 & 14000000 & 1142519 & 589629 & 2.7743e+003 & 2.7949e+003 & 130899 & 1.823e+000 & 1.077e+000 \\ 
32 & 12 & 0.6667 & 11000000 & 874736 & 453702 & 2.1478e+003 & 2.1269e+003 & 100289 & 1.838e+000 & 1.114e+000 \\ 
32 & 13 & 0.6667 & 14000000 & 1102885 & 571297 & 2.7163e+003 & 2.7051e+003 & 126588 & 1.846e+000 & 1.128e+000 \\ 
32 & 14 & 0.6667 & 12000000 & 938220 & 487280 & 2.2964e+003 & 2.2897e+003 & 107856 & 1.841e+000 & 1.131e+000 \\ 
32 & 15 & 0.6667 & 14000000 & 1086256 & 568979 & 2.6827e+003 & 2.6788e+003 & 125154 & 1.876e+000 & 1.159e+000 \\ 

32 & 1 & 1.0000 & 10000000 & 3332721 & 0      & 1.5819e+004 & 1.5825e+004 & 247677 & 1.719e+000 & 3.834e-001 \\ 
32 & 2 & 1.0000 & 10000000 & 2369493 & 0      & 1.1341e+004 & 1.1354e+004 & 179076 & 1.734e+000 & 5.347e-001 \\ 
32 & 2 & 0.5833 & 10000000 & 1381226 & 991383 & 1.9456e+003 & 1.8973e+003 & 178013 & 1.785e+000 & 5.536e-001 \\ 
32 & 2 & 0.5417 & 20000000 & 2573137 & 2181102 & 1.9085e+003 & 1.9266e+003 & 356466 & 1.771e+000 & 5.491e-001 \\ 
32 & 2 & 0.5208 & 40000000 & 4950753 & 4559821 & 1.9429e+003 & 1.8834e+003 & 712677 & 1.772e+000 & 5.496e-001 \\ 

64 & 1 & 0.6667 & 10000000 & 2233592 & 1123461 & 1.8555e+003 & 1.8491e+003 & 124756 & 1.708e+000 & 2.675e-001 \\ 
64 & 2 & 0.6667 & 20000000 & 3102155 & 1571597 & 2.6385e+003 & 2.6427e+003 & 175391 & 1.767e+000 & 3.935e-001 \\ 
64 & 4 & 0.6667 & 20000000 & 2208690 & 1127635 & 1.9124e+003 & 1.9001e+003 & 125776 & 1.806e+000 & 5.610e-001 \\ 
64 & 8 & 0.6667 & 22000000 & 1782288 & 917200 & 1.5424e+003 & 1.5407e+003 & 101851 & 1.825e+000 & 7.699e-001 \\ 
64 & 16 & 0.6667 & 30000000 & 1909428 & 988772 & 1.6775e+003 & 1.6584e+003 & 109624 & 1.855e+000 & 9.917e-001 \\ 
64 & 24 & 0.6667 & 25000000 & 1450673 & 755564 & 1.2598e+003 & 1.2479e+003 & 83414 & 1.847e+000 & 1.081e+000 \\ 

128 & 1 & 0.6667 & 10000000 & 2221748 & 1119251 & 6.5420e+002 & 6.5077e+002 & 62092 & 1.714e+000 & 1.906e-001 \\ 
128 & 16 & 0.6667 & 10000000 & 573255 & 297429 & 1.7767e+002 & 1.7490e+002 & 16466 & 1.851e+000 & 7.764e-001 \\ 

\\
\end{tabular}
\end{tiny}
\caption{Numerical results for $\{W_t^j\}$}
\label{table2}
\end{table}

\end{document}